\documentclass{article}
\usepackage{amsfonts, latexsym}

\newtheorem{theorem}{Theorem}

\begin{document}
\title{A Note on Cliques in Multipartite Graphs}
\author{
Julia V. Grishicheva\thanks{Maimonides State Classic Academy,
Moscow, Russia} , Alexander V. Seliverstov\thanks{Institute for
Information Transmission Problems, RAS, Moscow, Russia}\\
slvstv@iitp.ru }
\date{April 6, 2004}
\maketitle
\begin{abstract}
We consider a set of cliques in any multipartite graph with two
vertices in each part. Moreover, we construct a class of peculiar
polytopes.

Key words: multipartite graph, clique, polytope.
\end{abstract}

\paragraph{Introduction.} $n$-cliques in any $n$-partite graph
with two vertices in each part correspond to models for a
conjunctive normal form with two literals in each disjunction
(2CNF). Thus, the clique can be found in  polynomial time
\cite{EIS76}. On the other hand, some properties of a set of such
cliques are hard for  computing \cite{LS04}. Our article is a
contribution to this collection.

\paragraph{Terminology.}
A graph $\Gamma$ is said to be {\em multipartite} if the vertex
set is partitioned into {\em parts} such that there is no edge
between any pair of vertices from the same part. The multipartite
graph $\Gamma$ is {\em complete} if every two vertices from
different parts are joined. A complete multipartite graph $\Gamma$
is a complete graph if every part of $\Gamma$ consists of a single
vertex. A complete subgraph with exactly $n$ vertices is called a
$n$-{\em clique}.

There is a complete $3$-partite graph.
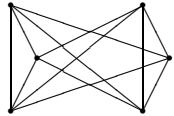
\begin{center}
\begin{picture}(60, 40)
\put(0, 0){\circle*{2}} \put(50, 0){\circle*{2}} \put(10,
20){\circle*{2}} \put(60, 20){\circle*{2}} \put(0,
40){\circle*{2}} \put(50, 40){\circle*{2}} \put(0, 0){\line(0,
1){40}} \put(0, 0){\line(1, 2){10}} \put(10, 20){\line(-1, 2){10}}
\put(50, 0){\line(0, 1){40}} \put(50, 0){\line(1, 2){10}} \put(60,
20){\line(-1, 2){10}} \put(0, 0){\line(3, 1){60}} \put(0,
0){\line(5, 4){50}} \put(10, 20){\line(2, -1){40}} \put(10,
20){\line(2, 1){40}} \put(0, 40){\line(3, -1){60}} \put(0,
40){\line(5, -4){50}}
\end{picture}
\end{center}

A closed and bounded subset $\Pi$ of the $d$-dimensional rational
affine space $\mathbb{Q}^d$ is called a {\em convex polytope} if
it is the set of solutions to a finite system of linear
inequalities. We shall omit convex for convex polytopes, and call
them simply polytopes. When the polytope is full-dimensional, each
nonredundant inequality corresponds to a facet.

Let $\Pi$ be a convex $d$-dimensional polytope in $\mathbb{Q}^d$.
To define {\em faces} geometrically, it is convenient to use
supporting hyperplanes. A hyperplane $H$ of $\mathbb{Q}^d$ is {\em
supporting} the polytope $\Pi$ if one of the closed half-spaces of
$H$ contains the polytope $\Pi$. A subset $\Phi \subseteq \Pi$ is
called a {\em face} of $\Pi$ if it is either empty, $\Pi$ itself
or the intersection of $\Pi$ with a supporting hyperplane.

The faces of dimension $0$, $1$, $d-2$ and $d-1$ are called the
{\em vertices, edges, ridges} and {\em facets}, respectively. The
vertices coincide with the extreme points of $\Pi$ which are
defined as points which cannot be represented as convex
combinations of two other points in the polytope $\Pi$. For a
subset $S \subseteq \mathbb{Q}^d$, the {\em convex hull} is
defined as the smallest convex set in $\mathbb{Q}^d$ containing
$S$. See also \cite{Sch86}.

\paragraph{Graphs and Polytopes.}
Let the indices $p$, $q$, $r$ are equal to either $1$ or $2$, and
the indices $i$, $j$, $k$ range over the segment $\{1, 2, \dots,
n\}$.

For any $n \ge 2$ we shall define a polytope $\Omega_n$ in the
$4n^2$ - dimensional rational affine space as a convex hull of
points $X$ with coordinates
$$ X_{ijpq} = \left \{
\begin{array} {cl}
1, & \mbox{if } p = \rho (i) \mbox{ and } q = \rho(j)\\
0, & \mbox{otherwise}.
\end{array} \right. $$
for any function $\rho: \{1, 2, \dots, n \}\to \{1, 2 \}$.

\begin{theorem}
Each vertex $X$ of the polytope $\Omega_n$ satisfy to the
equalities
\begin{eqnarray}
\forall i, j ~ \forall p, q ~ X_{ijpq} & = & X_{jiqp} \label{e1}\\
\forall i ~ X_{ii11} +X_{ii22} & = & 1 \label{e2} \\
\forall i ~ X_{ii12} & = & 0 \label{e3} \\
\forall i, j ~ \forall p ~ X_{ijp1} + X_{ijp2} & = &
X_{iipp}.\label{e4}
\end{eqnarray}
\end{theorem}
{\em Proof.} Any coordinate of a vertex $X$ is equal to either $0$
or $1$. Thus,
$$X_{ijpq} = 1 \mbox{ iff } X_{iipp} = 1 \mbox{ and } X_{jjqq} = 1;$$
$$
\mbox{if } i \ne j \mbox{ and } X_{ijpq} = 0, \mbox{ then }
X_{iipp} = 0 \mbox{ or } X_{jjqq} = 0.
$$

The equalities (\ref{e1}) as well as the equalities (\ref{e3}) are
obvious. The equalities (\ref{e2}) are satisfied because in each
sum one member is equal to $1$ and  others are equal to zero. The
equalities (\ref{e4}) are satisfied because in each sum no more
than one member is equal to $1$ and others are equal to zero.
$\Box$

With any point $X \in \Omega_n$ we associate a $n$-partite graph
with two vertices in each part, where for any pair $j \ne i$ of
indices the $p$-th vertex of the $i$-th part is joined by an edge
with the $q$-th vertex of the $j$-th part iff $X_{ijpq} > 0$. Any
positive coordinate of a point $X \in \Omega_n$ is convenient to
consider as the {\em weight} of a vertex or an edge.

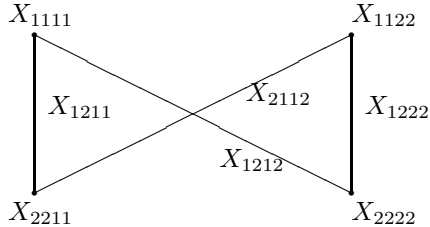
\begin{center}
\begin{picture}(140, 90)
\put(10, 10){\circle*{2}} \put(0, 0){$X_{2211}$} \put(130,
10){\circle*{2}} \put(130, 0){$X_{2222}$} \put(10,
70){\circle*{2}} \put(0, 75){$X_{1111}$} \put(130,
70){\circle*{2}} \put(130, 75){$X_{1122}$} \put(10, 10){\line(0,
1){60}} \put(10, 10){\line(2, 1){120}} \put(130, 10){\line(0,
1){60}} \put(130, 10){\line(-2, 1){120}} \put(15, 40){$X_{1211}$}
\put(135, 40){$X_{1222}$} \put(90, 45){$X_{2112}$} \put(80,
20){$X_{1212}$}
\end{picture}
\end{center}

Any vertex of the polytope $\Omega_n$ is a $n$-clique of the
complete $n$-partite graph with two vertices in each part.

\begin{theorem}
The dimension
$$\dim \Omega_n = \frac{n(n+1)}{2}.$$
\end{theorem}
{\em Proof.} The coordinates of a kind $X_{ij11}$ univocally
determine a point in the polytope $\Omega_n$. Thus, its dimension
$$\dim \Omega_n \le \frac{n(n+1)}{2}. $$
On the other hand, the polytope $\Omega_n$ contains a set of
$\frac{n(n+1)}{2}+1$ affine independent points. In more detail,
\begin{itemize}
\item [-] there is a point with every coordinate $X_{ii11}=0$;
\item [-] there are $n$ points with one unit and $n-1$ zeroes
among the coordinates $X_{ii11}$; \item [-] there are
$\frac{n(n-1)}{2}$ points with two units and $n-2$ zeroes among
the coordinates $X_{ii11}$.
\end{itemize}
Therefore, the dimension $\dim \Omega_n \ge \frac{n(n+1)}{2}$. As
a conclusion, we have the desired equality. $\Box$

\begin{theorem}
A convex hull of any vertex pair in the polytope $\Omega_n$ is an
edge.
\end{theorem}
{\em Proof.} Let us consider two vertices $X$ and $Y$ in the
polytope $\Omega_n$. They correspond to the pair of $n$-cliques
denoted $X$ and $Y$ too. Let us mark one edge in each clique that
does not belong to the other clique. Define a linear form
$$F(X) = \sum_{i>j} \alpha_{ijpq}X_{ijpq},$$
where $\alpha_{ijpq}$ is defined in this way. Let us consider an
edge between the $p$-th vertex in the $i$-th part and the $q$-th
vertex in the $j$-th part. Then
$$\alpha_{ijpq}=\left\{
\begin{array}{cl}
2, & \mbox{if this edge belongs to neither the $n$-clique $X$
nor the $n$-clique $Y$}\\
1, & \mbox{if this edge is marked}\\
0, & \mbox{otherwise}
\end{array}
\right.
$$
The equations $F(X)=1$ and $F(Y)=1$ are satisfied on these
vertices. And for any other vertex $Z$ we have $F(Z) \ge 2$. As a
consequence, these vertices $X$ and $Y$ are end points of an edge
in $\Omega_n$. $\Box$

\paragraph{Geometry of $\Omega_3$.}
Let us consider the polytope $\Omega_3$. The dimension $\dim
\Omega_3 = 6$. Any facet is determined by six of eight its
vertices, i.e. all except for two. These excluded vertices
correspond to a pair of triangles in a complete $3$-partite graph.
There are three cases.

{\em Excluded triangles have no common vertex.} Without loss of
generality we assume that the first triangle consists of the first
vertices and the second triangle consists of the second vertices
in each part.
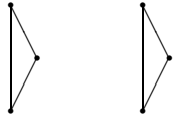
\begin{center}
\begin{picture}(60, 40)
\put(0, 0){\circle*{2}} \put(50, 0){\circle*{2}} \put(10,
20){\circle*{2}} \put(60, 20){\circle*{2}} \put(0,
40){\circle*{2}} \put(50, 40){\circle*{2}} \put(0, 0){\line(0,
1){40}} \put(0, 0){\line(1, 2){10}} \put(10, 20){\line(-1, 2){10}}
\put(50, 0){\line(0, 1){40}} \put(50, 0){\line(1, 2){10}} \put(60,
20){\line(-1, 2){10}}
\end{picture}
\end{center}
Other six vertices of the polytope $\Omega_3$ belong to a facet
defined by the equation
$$X_{1211}+X_{1311}+X_{2311}+X_{1222}+X_{1322}+X_{2322} = 1.$$

{\em Excluded triangles have a common edge.} Without loss of
generality we assume that first triangle consists of the first
vertices and the common edge joints the first part with the second
part.
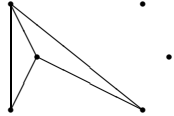
\begin{center}
\begin{picture}(60, 40)
\put(0, 0){\circle*{2}} \put(50, 0){\circle*{2}} \put(10,
20){\circle*{2}} \put(60, 20){\circle*{2}} \put(0,
40){\circle*{2}} \put(50, 40){\circle*{2}} \put(0, 0){\line(0,
1){40}} \put(0, 0){\line(1, 2){10}} \put(10, 20){\line(-1, 2){10}}
\put(50, 0){\line(-2, 1){40}} \put(50, 0){\line(-5, 4){50}}
\end{picture}
\end{center}
A convex hull of other six vertices of the polytope $\Omega_3$ is
a face defined by the equation $X_{1211} = 0$.

{\em Excluded triangles have a sole common vertex.} Without loss
of generality we assume that the first triangle consists of the
first vertices and the common vertex belongs to the first part.
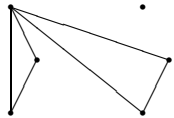
\begin{center}
\begin{picture}(60, 40)
\put(0, 0){\circle*{2}} \put(50, 0){\circle*{2}} \put(10,
20){\circle*{2}} \put(60, 20){\circle*{2}} \put(0,
40){\circle*{2}} \put(50, 40){\circle*{2}} \put(0, 0){\line(0,
1){40}} \put(0, 0){\line(1, 2){10}} \put(10, 20){\line(-1, 2){10}}
\put(50, 0){\line(1, 2){10}} \put(50, 0){\line(-5, 4){50}}
\put(60, 20){\line(-3, 1){60}}
\end{picture}
\end{center}
A convex hull of other six vertices of the polytope $\Omega_3$ can
not be a face. For example, the form
$$X_{1222} + X_{1312} + X_{1212} + X_{1221}$$
is evaluated as $0$ and $2$ on the excluded vertices of the
polytope $\Omega_3$, but it is evaluated as $1$ on other six
vertices.

\paragraph{Some remarks.}
The $3$-dimensional polytope $\Omega_2$ has four vertices, i.e.
$\Omega_2$ is a simplex.

Obviously, the polytope  $\Omega_n$ is symmetrical. Thus, each
vertex belongs to the same number of facets. The important
question is still open: how many facets are in the polytope
$\Omega_n$?

Note that in high dimensions, there are three {\em regular}
polytopes.
\begin{itemize}
\item [-] The $d$-dimensional simplex has $d+1$ vertices and $d+1$
facets. \item [-] The $d$-dimensional cube has $2^d$ vertices and
$2d$ facets. \item [-] The $d$-dimensional cross-polytope has $2d$
vertices and $2^d$ facets.
\end{itemize}
The {\em dual} of a regular polytope is another polytope, also
regular, having one vertex in the center of each facet of the
polytope we started with. The simplex is self-dual, and the dual
of the cube is the cross-polytope. So there is not a polynomial
upper bound for the number of facets intersecting at one vertex.

\paragraph{Acknowledgements.}
We wish to thank K.Yu.~Gorbunov for critical comments on this
article.

\end{document}